\documentclass[reqno]{slides}

\usepackage{amsmath}
\usepackage{amscd}
\usepackage{amssymb}
\usepackage{enumerate}
\usepackage{amsfonts}
\usepackage{graphicx}
\usepackage[all]{xy}
\usepackage{mathrsfs}

\usepackage{bbm}

\newcommand{\e}{\mathbbm{1}}

\newcommand{\rar}[1]{\stackrel{#1}{\longrightarrow}}
\newcommand{\xrar}[1]{\xrightarrow{#1}}

\newcommand{\into}{\hookrightarrow}

  \newcommand{\ga}{\gamma}
\newcommand{\Ga}{\Gamma}

\newcommand{\om}{\omega}

\newcommand{\bF}{{\mathbb F}}

\newcommand{\bK}{{\mathbb K}}

\newcommand{\bQ}{{\mathbb Q}}

\newcommand{\bZ}{{\mathbb Z}}

\newcommand{\cL}{{\mathcal L}}
\newcommand{\cM}{{\mathcal M}}

\newcommand{\cP}{{\mathcal P}}

\newcommand{\sC}{{\mathscr C}}
\newcommand{\sD}{{\mathscr D}}

\newcommand{\fg}{{\mathfrak g}}

\newcommand{\abs}[1]{\vert #1\vert}

\newcommand{\Hom}{\operatorname{Hom}}

\newcommand{\Ind}{\operatorname{Ind}}
\newcommand{\ind}{\operatorname{ind}}

\newcommand{\tens}{\otimes}

\newcommand{\Fun}{\operatorname{Fun}}

\newcommand{\Fr}{\operatorname{Fr}}

\newcommand{\ql}{\overline{\bQ}_\ell}

\newcommand{\Log}{\operatorname{Log}}

\begin{document}

\begin{slide}

\large

\centerline{\textbf{CHARACTER SHEAVES}}

\centerline{\textbf{ON UNIPOTENT GROUPS}}

\centerline{\textbf{IN CHARACTERISTIC} $p>0$}

\vskip2cm

\centerline{Mitya Boyarchenko}

\centerline{Vladimir Drinfeld}

\vskip1cm

\centerline{University of Chicago}

\end{slide}

\begin{slide}

\centerline{\small Overview}

{\tiny

These are slides for a talk given by the authors at the conference ``Current developments and directions in the Langlands program'' held in honor of Robert Langlands at the Northwestern University in May of 2008. The research program outlined in this talk was realized in a series of articles [1]--[4]. The orbit method for unipotent groups in positive characteristic is discussed in [1]. The results on character sheaves discussed in parts I and II of this talk are proved in [3]. The results described in part III are proved in [2] and [4]; the former studies $L$-packets of irreducible characters and the latter is devoted to the relationship between characters and character sheaves on unipotent groups over finite fields.

\begin{enumerate}[{[}1{]}]
\item M.~Boyarchenko and V.~Drinfeld, {\em A motivated introduction \\ to character sheaves and the orbit method for unipotent \\ groups in positive characteristic}, Preprint, {\tt math/0609769}
\item \vskip-0.3cm M.~Boyarchenko, {\em Characters of unipotent groups over \\ finite fields}, {\tt arXiv:0712.2614}, Selecta Math. \textbf{16} (2010), \\ no.~4, 857--933.
\item \vskip-0.3cm M.~Boyarchenko and V.~Drinfeld, {\em Character sheaves on \\ unipotent groups in positive characteristic: foundations}, \\ {\tt arXiv:0810.0794}, to appear in Selecta Mathematica.
\item \vskip-0.3cm M.~Boyarchenko, {\em Character sheaves and characters of \\ unipotent groups over finite fields}, {\tt arXiv:1006.2476}, \\ to appear in the American Journal of Mathematics.
\end{enumerate}

}

\end{slide}

\begin{slide}

\centerline{\large Some historical comments}

A geometric approach to representation theory for unipotent groups:
orbit method (Dixmier, Kirillov and others; late 1950s -- early
1960s). Does not work in characteristic $p>0$ unless additional
assumptions on the group are made.

The geometric approach to studying irreducible representations of
groups of the form $G(\bF_q)$, where $G$ is a \underline{reductive}
group over $\bF_q$, is at the heart of Deligne-Lusztig theory
(1970s) and Lusztig's theory of character sheaves (1980s).

In 2003 Lusztig explained that there should also exist an
interesting theory of character sheaves for unipotent groups in
char. $p>0$.

We will outline a theory that combines some of the essential
features of Lusztig's theory and of the orbit method.

\end{slide}

\begin{slide}

\centerline{\large\bf Part I. Definition of}
 \centerline{\large\bf character sheaves}

\vskip0.5cm

\centerline{Some notation}

$X=$ a scheme of finite type over a field $k$ \hfill\newline
 $G=$ an algebraic group over $k$ \hfill\newline
 $\ell\ =$ a prime different from $\operatorname{char}k$

$\sD(X):=D^b_c(X,\ql)$, the bounded derived \hfill\newline category
of constructible complexes of \hfill\newline $\ql$-sheaves on $X$

Given a $G$-action on $X$, one can define \\ $\sD_G(X)$, the
$G$-equivariant derived category.

In this talk we will mostly be concerned \\ with $\sD_G(G)$, defined
on the next slide, \\ where $G$ acts on itself by conjugation.


\end{slide}

\begin{slide}

\centerline{\large The category $\sD_G(G)$}

$G=$ a unipotent group over a field $k$ \hfill\newline $\ell=$ a
prime different from $\operatorname{char}k$

$\sD(G)=D^b_c(G,\ql)$, as before

$\mu:G\times G\rar{}G$ is the multiplication map \hfill\newline
$c:G\times G\rar{}G$ is the conjugation action \hfill\newline
$p_2:G\times G\rar{}G$ is the second projection

$\sD_G(G)$ consists of pairs $(M,\phi)$, where \hfill\newline
$M\in\sD(G)$ and $\phi:c^*M\rar{\simeq}p_2^*M$ satisfies \\ the
obvious cocycle condition

$\sD(G)$ is monoidal, and $\sD_G(G)$ is braided \hfill\newline
monoidal, with respect to \underline{convolution} \hfill\newline
(with compact supports):
\[
M*N=R\mu_!(M\boxtimes N)
\]
(Recall that a \emph{braiding} is a certain type \\ of a
commutativity constraint.)

\end{slide}

\begin{slide}

\centerline{\large Idempotents in monoidal categories}

(1) A \underline{weak idempotent} in a monoidal category
$(\cM,\otimes,\e)$ is an object $e$ such that $e\otimes e\cong e$.

(2) A \underline{closed idempotent} in $\cM$ is an object $e$ such
that there is an arrow $\pi:\e\rar{} e$, which becomes an
isomorphism after applying $e\otimes -$ as well as after applying
$-\otimes e$

It is more convenient to work with the second notion, because it is
much more rigid (e.g., $\pi$ is unique up to a unique automorphism
of $e$).

In each of the two contexts, we have the notion of the
\underline{Hecke subcategory} $e\cM e$. If $e$ is closed, $e\cM e$
is monoidal as well, with unit object $e$.

When $\cM$ is additive and braided, we can also talk about
\underline{minimal} (weak/closed) idempotents.

\end{slide}

\begin{slide}

\centerline{\large Motivation behind idempotents}

1. \textbf{Classical fact:} $\Ga=$ finite group; \hfill\newline
$\Fun(\Ga)=$ algebra of functions $\Ga\rar{}\ql$ w.r.t. pointwise
addition and convolution; \hfill\newline $\Fun(\Ga)^\Ga=$ center of
$\Fun(\Ga)$

$\exists$ a bijection between the set $\widehat{\Ga}$ of irreducible
characters of $\Ga$ over $\ql$ and the set of minimal idempotents in
$\Fun(\Ga)^\Ga$: namely,
$\chi\longmapsto\frac{\chi(1)}{\abs{\Ga}}\cdot\chi$

2. \textbf{Explanation of the term ``closed'':} Let $X$ be a scheme
of finite type over $k$, and let $\cM=\sD(X)$, with
$\otimes=\overset{L}{\otimes}_{\ql}$. Then the closed idempotents in
$\cM$ are of the form $i_!\ql$, where $i:Y\into X$ is the inclusion
of a \emph{closed} subscheme and $\ql$ is the constant sheaf on $Y$.

\emph{Typically (e.g., in this example) there are many weak
idempotents that are not closed.}

\end{slide}

\begin{slide}

3. \textbf{Orbit method:} $k=\overline{k}$,
$\operatorname{char}k=p>0$; \\ $G=$ connected unipotent group over
$k$; \\ assume $G$ has nilpotence class $<p$; \\
form $\fg=\Log G$ and $\fg^*=$ Serre dual of $\fg$

Then $G$ acts on $\fg$ and $\fg^*$, and we can consider
\[
\sD_G(G) \xrar{\exp^*} \sD_G(\fg) \xrar{\text{Fourier}}
\sD_G(\fg^*),
\]
which yields a bijection between minimal weak idempotents in
$\sD_G(G)$ and $G$-orbits in $\fg^*$.

Note that the $G$-orbits in $\fg^*$ are all \underline{closed} \\
(as $G$ is unipotent). Hence all minimal weak idempotents in
$\sD_G(G)$ are closed as well.

In fact, the last statement remains true for \underline{any}
unipotent group $G$ over $k$.

\end{slide}

\begin{slide}

\centerline{\large Character sheaves and $L$-packets}

$k=$ an algebraically closed field of char. $p>0$ \\
$G=$ an arbitrary unipotent group over $k$

Pick a minimal closed idempotent $e\in\sD_G(G)$. Note that
$e\sD_G(G)e=e\sD_G(G)$. Consider
\[
\cM_e^{perv}=\{(M,\phi)\in e\sD_G(G)\,\lvert\, M\text{ is
perverse}\},
\]
a full additive subcategory of $\sD_G(G)$.

\textbf{Definition.} The \underline{$L$-packet of character sheaves}
associated to $e$ is the set of indecomposable objects of
$\cM_e^{perv}$. An object of $\sD_G(G)$ is a \underline{character
sheaf} if it lies in some $L$-packet.

\textbf{Conjecture.} (1) $\cM^{perv}_e$ is a semisimple abelian
category with finitely many simple objects. \\
(2) $\exists n_e\in\bZ$ such that $e[-n_e]\in\cM_e^{perv}$. \\
(3) $\cM_e:=\cM_e^{perv}[n_e]$ is closed under $*$, and is a
\emph{modular} ($\approx$as far from being symmetric as possible)
braided monoidal category.

\end{slide}

\begin{slide}

\centerline{\large\bf Part II. Construction}
 \centerline{\large\bf of character sheaves}

\vskip0.5cm

\centerline{\large Overview}

\textbf{Two classical results.} (1) If $\Ga$ is a finite nilpotent
group and $\rho\in\widehat{\Ga}$, there exist a subgroup
$H\subset\Ga$ and a homomorphism $\chi:H\to\ql^\times$ such that
$\rho=\rho_{H,\chi}:=\Ind_H^\Ga\chi$.

(2) For \emph{any} such pair $(H,\chi)$, \underline{Mackey's
criterion} states that $\rho_{H,\chi}$ is irreducible $\iff$
$\forall\,\ga\in\Ga\setminus H$,
\[
\chi\bigl\lvert_{H\cap \ga H\ga^{-1}} \not\equiv
\chi^{\ga}\bigl\lvert_{H\cap \ga H\ga^{-1}}
\]

The notion of a $1$-dimensional representation, the operation of
induction, and the two results stated above, have geometric
analogues. This is what the second part will be about.

\end{slide}

\begin{slide}

\centerline{\large $1$-dimensional character sheaves}

A geometrization of the notion of a $1$-dimensional representation
is provided by the notion of a multiplicative local system.

$G=$ an algebraic group over a field $k$ \\
$\ell=$ a prime different from $\operatorname{char}k$ \\
$\mu:G\times G\rar{} G$ is the multiplication map

A nonzero $\ql$-local system $\cL$ on $G$ is said to be
\underline{multiplicative} if $\mu^*\cL\cong\cL\boxtimes\cL$ (so
$\operatorname{rk}\cL=1$).

If $G$ is connected and unipotent, and $\bK_G$ is the dualizing
complex of $G$, then $e_{\cL}:=\cL\otimes\bK_G$ is a minimal closed
idempotent in $\sD_G(G)$.

If, moreover, $k=\overline{k}$, the corresponding $L$-packet
consists of the single character sheaf $\cL[\dim G]$.

\end{slide}

\begin{slide}

\centerline{\large A more canonical viewpoint}

Fix an embedding $\psi:\bQ_p/\bZ_p\into\ql^\times$. If $\Ga$ is a
finite $p$-group, every homomorphism $\Ga\to\ql^\times$ factors
through $\psi$. Note that $\Hom(\Ga,\bQ_p/\bZ_p)$ does not depend on
$\ell$.

Now let $G$ be an algebraic group over a field $k$ of char. $p>0$.
We have the functors

$\bigl\{$central extensions $1\to\bQ_p/\bZ_p\to\widetilde{G}\to G\to
1\bigr\}$ \\ $\xrightarrow{forgetful}\bigl\{\text{multiplicative }
\bQ_p/\bZ_p-\text{torsors on } G\bigr\}$ \\ $\xrightarrow{\ \
\psi_*\ \ } \bigl\{\text{multiplicative }\ql-\text{local systems on
} G\bigr\}$

If $G$ is \emph{connected} and \emph{unipotent}, they induce
bijections on isomorphism classes of objects.

So we will study central extensions by $\bQ_p/\bZ_p$ in place of
multiplicative local systems. \\ The composition of the two functors
will be \\ denoted by $\chi\longmapsto\cL_\chi$.

\end{slide}

\begin{slide}

\centerline{\large Serre duality}

$k=$ perfect field of char. $p>0$ \\
$G=$ connected unipotent group over $k$

The \emph{Serre dual} of $G$ is the functor
\[
G^* : \left\{
\begin{array}{c}
\text{perfect} \\
k-\text{schemes}
\end{array}
\right\} \rar{} \left\{
\begin{array}{c}
\text{abelian} \\
\text{groups}
\end{array}
\right\}
\]
\[
S\longmapsto \bigl\{ \text{iso. classes of central extensions}
\]
\[
 \text{of the group scheme } G\times_k S \text{ by } \bQ_p/\bZ_p
\times S \bigr\}
\]

\textbf{Classical Serre duality.} If $G$ is commutative, then $G^*$
is representable by a perfect connected commutative unipotent group
scheme over $k$.

\textbf{Proposition.} In general, $G^*$ is representable by a
possibly disconnected perfect commutative unipotent group over $k$.
Its neutral connected component, $(G^*)^\circ$, can be naturally
identified with $(G^{ab})^*$.

\end{slide}

\begin{slide}

\centerline{\large Definition of admissible pairs}

$k=$ algebraically closed field of char. $p>0$ \\
$G=$ connected unipotent group over $k$

An \underline{admissible pair} for $G$ is a pair $(H,\chi)$, where
$H\subset G$ is a connected subgroup and $\chi\in H^*(k)$, such that
the following three conditions hold:

(1) Let $G'$ be the normalizer in $G$ of the pair $(H,\chi)$; then
$G^{\prime\circ}/H$ is commutative.

(2) The homomorphism $G^{\prime\circ}/H\to (G^{\prime\circ}/H)^*$
induced by $\chi$ (a geometrization of $g\mapsto\chi([g,-])$) is an
isogeny (i.e., has finite kernel).

(3) If $g\in G(k)\setminus G'(k)$, the restrictions of $\chi$ and
$\chi^g$ to $(H\cap g^{-1}Hg)^\circ$ are nonisomorphic.

This is the correct geometric analogue of \hfill\newline Mackey's
irreducibility criterion.

\end{slide}

\begin{slide}


\centerline{\large Induction with compact supports}

$G=$ unipotent group over a field $k$ \\
$G'\subset G$ is a closed subgroup

One can define a functor
\[
\ind_{G'}^G : \sD_{G'}(G')\rar{}\sD_G(G),
\]
called \underline{induction with compact supports}.

Its construction is standard: take $M\in\sD_{G'}(G')$, extend by
zero to $\overline{M}\in\sD(G)$, then average, \emph{in the sense of
``lower shriek''}, with respect to the conjugation action of $G$ on
itself.

If $k$ is finite \emph{and $G'$ is connected}, this functor is
compatible with induction of functions via the sheaves-to-functions
dictionary. (If $G'$ is not connected, this may be false, in
general.)



\end{slide}

\begin{slide}

\centerline{\large Construction of $L$-packets}

$k=$ algebraically closed field of char. $p>0$ \\
$G=$ connected unipotent group over $k$ \\
$(H,\chi)=$ an admissible pair for $G$ \\
$G'=$ normalizer of $(H,\chi)$ in $G$

\underline{Define}: $\cL_\chi=$ multiplicative $\ql$-local system on
$H$ arising from $\chi$ via a chosen
$\psi:\bQ_p/\bZ_p\into\ql^\times$;

\vskip-0.7cm

$e_\chi = \cL_\chi\otimes\bK_H \cong\cL_\chi[2\dim H] \in
\sD_{G'}(H)$;

\vskip-0.7cm

$\sD_{G'}(G')\ni e = $ extension of $e_\chi$ by zero;

\vskip-0.7cm

$\sD_G(G)\ni f=\ind_{G'}^G e$

\textbf{Theorem.} (1) $f$ is a minimal closed \\ idempotent in
$\sD_G(G)$.

\vskip-0.7cm

(2) $\ind_{G'}^G$ restricts to a braided monoidal \\
equivalence
$e\sD_{G'}(G')\rar{\sim}f\sD_G(G)$. \\



\end{slide}

\begin{slide}

(3) Set $n_e=\dim H$, $n_f=\dim H-\dim(G/G')$,
\[
\cM_e=\cM_e^{perv}[n_e], \quad \cM_f=\cM_f^{perv}[n_f].
\]
Then $\cM_e$ and $\cM_f$ are closed under $*$, and
\[\ind_{G'}^G(\cM_e)=\cM_f.\]

(4) The categories $\cM_e$ and $\cM_f$ are semisimple and have
finitely many simple objects.

It follows that the $L$-packet of character sheaves on $G$
associated to $f$ is formed by the objects $\ind_{G'}^G(M_j)[-n_f]$,
where $M_1,\dotsc,M_k\in\cM_e$ are the simple objects. The objects
$M_j$ can be described explicitly.

(5) Every minimal \emph{weak} idempotent $f\in\sD_G(G)$ arises from
an admissible pair $(H,\chi)$ as above, and, in particular, is
necessarily \emph{closed}.

\emph{However, $(H,\chi)$ is usually far from unique.}

\end{slide}

\begin{slide}

\centerline{\large\bf Part III. Relation to characters}

\vskip0.5cm

\centerline{\large $L$-packets of irreducible characters}

$\bF_q=$ finite field with $q$ elements \\
$G_0=$ connected unipotent group over $\bF_q$ \\
$G=G_0\otimes_{\bF_q}\overline{\bF}_q$ \\
$\Fr : G\to G$ is the Frobenius endomorphism

$\cP=$ a $\Fr$-stable $L$-packet of character sheaves on $G$
(equivalently, the corresponding minimal idempotent $e\in\sD_G(G)$
satisfies $\Fr^*e\cong e$)

Define
\[
\cP^{\Fr}=\bigl\{ N\in\cP \,\bigl\lvert\, \Fr^*(N)\cong N\bigr\}.
\]
For each $M\in\cP^{\Fr}$, form the corresponding function
$t_M:G_0(\bF_q)\rar{}\ql$ (it is well defined up to rescaling).

\end{slide}

\begin{slide}

We define a set $\cP'$ of irreducible characters of $G_0(\bF_q)$ as
follows: $\om\in\cP'$ $\iff$ $\om$ lies in the span of the set of
functions $\bigl\{t_M\bigr\}_{M\in\cP^{\Fr}}$.

We call $\cP'\subset\widehat{G_0(\bF_q)}$ the $L$-packet of
irreducible characters defined by the $\Fr$-stable $L$-packet $\cP$
of character sheaves.

\textbf{Theorem.} Every irreducible character of $G_0(\bF_q)$ over
$\ql$ lies in an $L$-packet as defined above.

There exists a description of $L$-packets of \\ irreducible
characters of $G_0(\bF_q)$ in terms of \\ admissible pairs for
$G_0$, which is analogous \\ to the statement we discussed earlier,
but is \\ independent of the theory of character sheaves.

This description plays an important role in \\ the proof of the last
theorem.

\end{slide}

\begin{slide}

\centerline{\large Description of $L$-packets of characters}

Let $G_0$ be a connected unipotent group over $\bF_q$, and consider
all pairs $(H,\chi)$ consisting of a connected subgroup $H\subset
G_0$ and an element $\chi\in H^*(\bF_q)$ such that
$(H\tens_{\bF_q}\overline{\bF}_q,\chi\tens_{\bF_q}\overline{\bF}_q)$
is an admissible pair for $G=G_0\tens_{\bF_q}\overline{\bF}_q$.

Two such pairs, $(H_1,\chi_1)$ and $(H_2,\chi_2)$, are
\hfill\newline \underline{geometrically conjugate} if they are
conjugate by an element of $G(\overline{\bF}_q)$.


Let $\sC$ be a geometric conjugacy class of such pairs. For an irrep
$\rho$ of $G_0(\bF_q)$ over $\ql$, we will write $\rho\in \cP_\sC$
if there exists $(H,\chi)\in\sC$ such that $\rho$ is an irreducible
constituent of $\Ind_{H(\bF_q)}^{G_0(\bF_q)}t_{\cL_\chi}$.

\textbf{Theorem.} The $\cP_\sC$ are exactly the $L$-packets of
irreducible characters of $G_0(\bF_q)$.

It could happen that $\sC_1\neq\sC_2$ and $\cP_{\sC_1}=\cP_{\sC_2}$.

\end{slide}

\end{document}